\documentclass[12pt]{article}
\usepackage{amssymb,amsthm,amsmath,latexsym}
\usepackage{url}
\usepackage{color}
\usepackage{comment}
\usepackage{lscape}
\usepackage{mathtools}
\usepackage{enumitem}
\usepackage{graphicx}
\newtheorem{thm}{Theorem}

\newtheorem{lem}[thm]{Lemma}

\theoremstyle{remark}
\newtheorem{rem}[thm]{Remark}
\newcommand{\FF}{\mathbb{F}}
\newcommand{\ZZ}{\mathbb{Z}}

\newcommand{\cD}{\mathcal{D}}
\newcommand{\cP}{\mathcal{P}}
\newcommand{\cB}{\mathcal{B}}

\newcommand{\0}{\mathbf{0}}
\newcommand{\1}{\mathbf{1}}
\renewcommand{\vec}[1]{\overrightarrow{#1}}
\DeclareMathOperator{\wt}{wt}

\newcommand{\binpart}[1]{\cB(#1)}
\allowdisplaybreaks

\begin{document}

\title{Nonexistence of a quasi-symmetric $2$-$(33,15,35)$ design}

\author{
Masaaki Harada\thanks{
Research Center for Pure and Applied Mathematics,
Graduate School of Information Sciences,
Tohoku University, Sendai 980--8579, Japan.
email: \protect\url{mharada@tohoku.ac.jp}
}
}

\maketitle

\begin{abstract}
It is shown that there exists no quasi-symmetric $2$-$(33,15,35)$ design
with  intersection numbers $6$ and $9$.
This is proved by considering
the ternary code generated by the incidence vectors
of the blocks and its dual code.
\end{abstract}

\section{Introduction}

A \emph{$2$-$(v,k,\lambda)$ design} $\cD$ is a pair consisting of a set $\cP$ of $v$ points
and a collection $\cB$ of $k$-element subsets of 
$\cP$ (called blocks)
such that every $2$-element subset of $\cP$ is contained in exactly
$\lambda$ blocks.
The number of blocks that contain a given
point is traditionally denoted by $r$, and the total number of
blocks is $b$.
We assume that $2 < k < v$ throughout this paper.
A design with no repeated block is called \emph{simple}.
In this paper, we consider only simple designs.
Often a $2$-$(v,k,\lambda)$ design is simply called a $2$-design.
A $2$-design can be represented by its (block by point)
\emph{incidence matrix} $A=(a_{ij})$, where $a_{ij}=1$ if the $j$th
point is contained in the $i$th block
and $a_{ij}=0$ otherwise.

A $2$-design is called \emph{quasi-symmetric} if 
the number
of points in the intersection of any two distinct blocks 
takes only two values $x$ and $y$ $(x < y)$.
The values $x$ and $y$ are called the \emph{intersection numbers}.
Quasi-symmetric $2$-designs are closely related to
strongly regular graphs and codes (see, e.g., \cite{BVM}, \cite{SS} and \cite{T07}).
Four classes of quasi-symmetric $2$-designs were given 
in~\cite{N}, namely, multiples of symmetric $2$-designs,
Steiner systems with $v > k^2$, strongly resolvable $2$-designs
and residual designs of biplanes.
A quasi-symmetric $2$-design $\cD$ is called \emph{exceptional} if 
neither $\cD$ nor its complementary design belongs to one of these four classes.
It is a fundamental problem to determine whether there exists
an exceptional quasi-symmetric $2$-design for a given set of parameters.
The feasible parameters of exceptional quasi-symmetric $2$-$(v,k,\lambda)$ designs 
with $2k \le v \le 70$
are listed in~\cite[Table~48.25]{QSHandbook}.
Codes have been used to investigate these existence problems
 (see, e.g., \cite{GPD}, \cite{HMT}, \cite{MT} and \cite{MT2}).
For recent nonexistence results based on restrictions
on block graphs, see~\cite{PSR} and \cite{RP}.

The existence of a quasi-symmetric $2$-$(33,15,35)$ design
is listed as an open question in~\cite[Table~48.25]{QSHandbook}
(see also~\cite[Section~8.5.6]{BVM}).
Such a design necessarily has intersection numbers $6$ and $9$
(see~\cite[Section~8.5.6]{BVM}).
Among the parameter sets for exceptional quasi-symmetric designs
listed as unresolved in that reference,
this one has the smallest number of points.
The aim of this paper is to establish the following theorem.

\begin{thm}\label{main}
There exists no quasi-symmetric $2$-$(33,15,35)$ design.
\end{thm}

\begin{rem}
If such a design existed, its block graph, in which two distinct
blocks are adjacent if and only if they intersect in exactly  $9$
points, would be a strongly regular graph with parameters
$(176,45,18,9)$.
Strongly regular graphs with these parameters are known to exist (see~\cite{AEB}).
\end{rem}

\section{Preliminaries}

Our approach is to consider the ternary code generated by
the incidence vectors of the blocks of a hypothetical design,
together with its dual code.
In this section, we introduce the notation and definitions
used in the proof.
Throughout the paper, the ring or field in which a calculation
is performed is understood from the context.
We specify it explicitly whenever necessary to avoid ambiguity.

\subsection{Ternary codes}
Let $\FF_3=\{0,1,2\}$ denote the finite field of order $3$.
A \emph{ternary code} $C$ of \emph{length} $n$
is a vector subspace of $\FF_3^n$.
A vector of $C$ is called a \emph{codeword} of $C$.
The \emph{weight} $\wt(x)$ of a vector $x \in \FF_3^n$ is
the number of nonzero components of $x$.

The \emph{dual} code $C^{\perp}$ of a ternary code $C$ of length $n$
is defined as follows:
\[
C^{\perp}=
\{x \in \FF_3^n \mid  x \cdot y = 0 \text{ for all } y \in C\},
\]
where
\[
x \cdot y = \sum_{i=1}^{n} x_i {y_i}
\]
for $x=(x_1,x_2,\ldots,x_n)$ and $y=(y_1,y_2,\ldots,y_n) \in \FF_3^n$.
A ternary code $C$ is called \emph{self-orthogonal} if $C \subseteq C^\perp$.

\subsection{Maps and inner product}

We employ the following maps and inner product.
Define a map $\phi:\FF_3\to\ZZ$ by
\[
 \phi(0)=0, \phi(1)=1 \text{ and }\phi(2)=-1.
\]
Let $\rho:\ZZ\to\FF_3$ denote the reduction map modulo $3$.
Then
\[
 \rho(\phi(a))=a
\]
for every element $a\in\FF_3$.
We extend both maps coordinatewise to vectors
and use the same symbols for these extensions.
Thus, for $u=(u_1,u_2,\ldots,u_n)\in\FF_3^n$ and 
for $z=(z_1,z_2,\ldots,z_n)\in\ZZ^n$, we write
\begin{align*}
 \phi(u)&=(\phi(u_1),\phi(u_2),\ldots,\phi(u_n))\in\ZZ^n \text{ and}\\
 \rho(z)&=(\rho(z_1),\rho(z_2),\ldots,\rho(z_n))\in\FF_3^n.
\end{align*}
In particular, $\rho(\phi(u))=u$ for every vector $u\in\FF_3^n$.
For vectors $a=(a_1,a_2,\ldots,a_n)$ and
$b=(b_1,b_2,\ldots,b_n) \in \ZZ^n$,
define their integer inner product by
\[
 \langle a,b\rangle_{\ZZ}=\sum_{i=1}^{n}a_i b_i.
\]

\section{Proof of Theorem~\ref{main}}\label{Sec:proof}

The aim of this section is to give a proof of Theorem~\ref{main}
that is as elementary as possible.

\subsection{Ternary code associated with the design}

Suppose that there exists a quasi-symmetric $2$-$(33,15,35)$ design
$\cD=(X,\cB)$
with  intersection numbers $6$ and $9$.
We may assume that $X=\{1,2,\ldots,33\}$.
For each subset $W\subseteq X$,
let $\vec{W}\in\{0,1\}^{33}$ denote its characteristic vector.

Let $A$ be the $176\times33$ incidence matrix of $\cD$.
Then the $(B,p)$-entry of $A$ is $1$ if $p\in B$ and $0$ otherwise,
and the row indexed by $B$ is $\vec{B}$.
Since $r=80$ and  $\lambda=35$, we have
\begin{equation}\label{eq:incidence}
 A^T A=45I_{33}+35J_{33} \text{ and } 
 \sum_{B\in\cB}\vec{B}=80\1
\end{equation}
over $\ZZ$, where $A^T$ denotes the transpose of $A$,
$I_{33}$ is the identity matrix of order $33$,
$J_{33}$ is the $33\times33$ all-ones matrix
and  $\1$ denotes the all-ones vector of length $33$.

Let $C(\cD)$ be the ternary code generated by the rows of $A$, that is,
\[
 C(\cD)=\langle \vec{B} \mid B\in\cB\rangle
 \subseteq\FF_3^{33}.
\]
The $176\times176$ integer matrix $AA^T$ has diagonal entries $15$
and off-diagonal entries $6$ or $9$, since $\cD$
has intersection numbers $6$ and $9$.
Since all these entries are divisible by $3$, we have
$AA^T\equiv O \pmod3$, where $O$ is the $176\times176$ zero matrix.
Thus, $C(\cD)$ is self-orthogonal.
Reducing the second equality in~\eqref{eq:incidence} modulo $3$,
we obtain
\[
C(\cD) \ni
\sum_{B\in\cB}\vec{B}=2\1.
\]
It follows that
\begin{equation}\label{eq:one}
 \1\in C(\cD) \subseteq C(\cD)^\perp.
\end{equation}

\subsection{Restrictions on $(0,1)$-vectors in $C(\cD)^\perp$}

Let $u=(u_1,u_2,\ldots,u_{33})\in C(\cD)^\perp$.
For each block $B\in\cB$, let
\[
 t(u)=\sum_{p=1}^{33}\phi(u_p)\in\ZZ
\text{ and }
 S_B=\sum_{p\in B}\phi(u_p)\in\ZZ.
\]

\begin{lem}\label{lem:moments}
Let $u$ be a vector of $C(\cD)^\perp$.
Then
\begin{equation}\label{eq:moments}
\begin{aligned}
 S_B &\in3\ZZ,\\
 \sum_{B\in\cB}S_B&=80t(u) \text{ and }\\
 \sum_{B\in\cB}S_B^2&=45\wt(u)+35t(u)^2.
\end{aligned}
\end{equation}
\end{lem}

\begin{proof}
Since $\rho$ is additive and $\rho(\phi(a))=a$
for every element $a\in\FF_3$, we have
\[
 \rho(S_B)
 =\sum_{p\in B}u_p
 =\sum_{p=1}^{33}u_p(\vec{B})_p
 =u\cdot\vec{B}=0
\]
for each block $B\in\cB$.
Here, $(\vec{B})_p$ denotes the $p$th coordinate of $\vec{B}$.
The last equality follows from
$u\in C(\cD)^\perp$ and $\vec{B}\in C(\cD)$.
Thus, $S_B$ is divisible by $3$ for every block $B\in\cB$.

Since $r=80$, we have
\[
 \sum_{B\in\cB}S_B
 =\sum_{p=1}^{33}\phi(u_p)\,|\{B\in\cB \mid p\in B\}|
 =80t(u).
\]
The entry of $A\phi(u)^T$ indexed by $B$ is $S_B$.
From~\eqref{eq:incidence}, we obtain
\begin{align*}
 \sum_{B\in\cB}S_B^2
 &=(A\phi(u)^T)^T(A\phi(u)^T)\\
 &=\phi(u)A^TA\phi(u)^T\\
 &=45\sum_{p=1}^{33}\phi(u_p)^2
   +35\left(\sum_{p=1}^{33}\phi(u_p)\right)^2\\
 &=45\wt(u)+35t(u)^2.
\end{align*}
The last equality follows from
$\phi(0)=0$ and $\phi(1)^2=\phi(2)^2=1$.
This completes the proof.
\end{proof}


For use in the following lemma, note that
\[
 S_B=\sum_{p\in B}\phi(u_p)
     =\langle\phi(u),\phi(\vec{B})\rangle_{\ZZ}
\]
for every vector $u\in C(\cD)^\perp$ and $B\in\cB$.

\begin{lem}\label{lem:signed}
Suppose that $u\in C(\cD)^\perp$ and that $\phi(u)$ has
$6$ entries equal to $1$, $3$ entries equal to $-1$ and $24$ entries equal to $0$.
Then 
\[
 S_B\in\{0,3\}
\]
for every block $B\in\cB$.
\end{lem}

\begin{proof}
Since $\phi(u)$ has $6$ entries equal to $1$,
$3$ entries equal to $-1$ and all remaining entries equal to $0$, and
$S_B=\langle\phi(u),\phi(\vec{B})\rangle_{\ZZ}$,
we have $-3\le S_B\le6$.
By Lemma~\ref{lem:moments}, we have $S_B\in3\ZZ$.
It follows that
\[
 S_B\in\{-3,0,3,6\}.
\]
Thus, for every  block $B\in\cB$, we have
\[
 S_B(S_B-3)\ge0,
\]
with equality if and only if $S_B\in\{0,3\}$.

Moreover, 
since $\wt(u)=9$ and $t(u)=3$, Lemma~\ref{lem:moments} gives
\begin{align*}
 \sum_{B\in\cB}S_B&=240 \text{ and }\\
 \sum_{B\in\cB}S_B^2&=45\cdot9+35\cdot3^2=720.
\end{align*}
Hence, we have
\[
 \sum_{B\in\cB}S_B(S_B-3)=720-3\cdot240=0.
\]
Therefore, we conclude that $S_B\in\{0,3\}$ for every block $B\in\cB$.
\end{proof}

We now apply Lemma~\ref{lem:moments} to characteristic vectors
of subsets of $X$.
Suppose that $W\subseteq X$ satisfies $\vec{W}\in C(\cD)^\perp$.
Let $u=\vec{W}$ and $w=|W|$.
The vector $u$ is a $(0,1)$-vector of weight $w$.
Since $\phi(0)=0$ and $\phi(1)=1$, we have
\[
 t(u)=\sum_{p=1}^{33}\phi(u_p)=w.
\]
Moreover, we have
\[
 S_B=\sum_{p\in B}\phi(u_p)=|B\cap W|
\]
for each block $B\in\cB$.
By Lemma~\ref{lem:moments}, $S_B=|B\cap W|$ is divisible by $3$.
Since $|B|=15$, we also have
\[
 0\le |B\cap W|\le\min(w,15).
\]
Let
\[
 \mathcal{I}_w=\{i\in3\ZZ\mid 0\le i\le\min(w,15)\}.
\]
Thus, $|B\cap W|\in\mathcal{I}_w$ for every $B\in\cB$.
For each $i\in\mathcal{I}_w$, define
\begin{equation}\label{eq:binarymoments}
 n_i=|\{B\in\cB \mid |B\cap W|=i\}|.
\end{equation}
Thus, $n_i$ is the number of blocks meeting $W$ in exactly $i$ points.

Since
\[
n_i
=|\{B\in\cB\mid |B\cap W|=i\}|
=|\{B\in\cB\mid S_B=i\}|,
\]
each block is counted exactly once in $\sum_{i\in\mathcal{I}_w} n_i$.
Thus, we have
\begin{equation}\label{eq:distribution0}
 \sum_{i\in\mathcal{I}_w} n_i=|\cB|=176.
\end{equation}
For each $i\in\mathcal{I}_w$, the $n_i$ blocks with $S_B=i$
contribute $i n_i$ to $\sum_{B\in\cB}S_B$
and $i^2 n_i$ to $\sum_{B\in\cB}S_B^2$.
Hence, we obtain
\begin{equation}\label{eq:distribution1}
\begin{aligned}
 \sum_{B\in\cB}S_B
 &=\sum_{i\in\mathcal{I}_w}i n_i \text{ and }\\
 \sum_{B\in\cB}S_B^2
 &=\sum_{i\in\mathcal{I}_w}
   \sum_{\substack{B\in\cB\\S_B=i}}S_B^2
 =\sum_{i\in\mathcal{I}_w}i^2 n_i.
\end{aligned}
\end{equation}
Using~\eqref{eq:distribution0} and~\eqref{eq:distribution1},
and applying Lemma~\ref{lem:moments} with $\wt(u)=t(u)=w$,
we obtain
\begin{equation}\label{eq:distribution}
\begin{aligned}
 \sum_{i\in\mathcal{I}_w} n_i &=176,\\
 \sum_{i\in\mathcal{I}_w} i n_i &=80w \text{ and }\\
 \sum_{i\in\mathcal{I}_w} i^2 n_i &=45w+35w^2,
\end{aligned}
\end{equation}
where all equalities are taken over $\ZZ$.

The following lemma restricts the possible weights of
$(0,1)$-vectors in $C(\cD)^\perp$.

\begin{lem}\label{lem:weights}
Suppose that $x\in C(\cD)^\perp\cap\{0,1\}^{33}$.
Then
\begin{equation}\label{eq:weights}
 \wt(x)\in\{0,15,18,33\}.
\end{equation}
\end{lem}

\begin{proof}
Since $\1\in C(\cD)$ by~\eqref{eq:one}
and $x\in C(\cD)^\perp$, we have
\[
 \rho(\wt(x))=x\cdot\1=0.
\]
Thus, $\wt(x)$ is divisible by $3$.
Since $\1\in C(\cD)\subseteq C(\cD)^\perp$ by~\eqref{eq:one}
and $x\in C(\cD)^\perp$,
the linearity of $C(\cD)^\perp$ gives
$\1+2x\in C(\cD)^\perp$.
Moreover, since $x$ is a $(0,1)$-vector, we have
\[
 \1+2x\in\{0,1\}^{33}
 \text{ and }
 \wt(\1+2x)=33-\wt(x).
\]
Therefore, it suffices to rule out weights $\wt(x)=3$, $6$, $9$ and $12$.
We consider these weights in the order $\wt(x)=3$, $9$, $6$ and $12$,
since the later cases use the exclusions established in the earlier ones.

\begin{enumerate}
\item
Suppose that $\wt(x)=3$.  Take $u=x$ in
Lemma~\ref{lem:moments}.
Then, since $x$ is a $(0,1)$-vector, we have
$\wt(u)=t(u)=3$ and $0\le S_B\le3$ for every block $B\in\cB$.
Since $S_B$ is divisible by $3$, we have $S_B\in\{0,3\}$.
In either case, $S_B^2=3S_B$.

On the other hand, \eqref{eq:moments} gives
\begin{align*}
 \sum_{B\in\cB}S_B&=80\cdot3=240 \text{ and }\\
 \sum_{B\in\cB}S_B^2&=45\cdot3+35\cdot3^2=450.
\end{align*}
This contradicts
\[
 \sum_{B\in\cB}S_B^2
 =\sum_{B\in\cB}3S_B
 =3\sum_{B\in\cB}S_B
 =720.
\]

\item
Suppose that $\wt(x)=9$.  Take $u=x$ in
Lemma~\ref{lem:moments}.
Then $\wt(u)=t(u)=9$ and $0\le S_B\le9$ for every block $B\in\cB$.
Since $S_B$ is divisible by $3$, we have
\[
 S_B\in\{0,3,6,9\}.
\]
Then we have 
\[
(S_B-3)(S_B-6)\ge0.
\]

On the other hand, \eqref{eq:moments} and $b=176$ give
\begin{align*}
\sum_{B\in\cB}(S_B-3)(S_B-6)
 &=\sum_{B\in\cB}S_B^2
   -9\sum_{B\in\cB}S_B+18b\\
 &=\bigl(45\cdot9+35\cdot9^2\bigr)
   -9(80\cdot9)+18\cdot176\\
 &=-72.
\end{align*}
This is a contradiction.

\item
Suppose that $\wt(x)=6$.
Write $x=\vec{W}$, where $W\subseteq X$ and $|W|=6$.
Let $n_0,n_3$ and $n_6$ be as defined in~\eqref{eq:binarymoments}.
From~\eqref{eq:distribution}, we have
\begin{align*}
 n_0+n_3+n_6&=176,\\
 3n_3+6n_6&=480 \text{ and }\\
 9n_3+36n_6&=1530.
\end{align*}
This  gives $18n_6=90$.
Since $n_6>0$, there is a block $B\in\cB$ with
$|B\cap W|=6$.  Hence, we have $W\subseteq B$.
By the definition of $C(\cD)$ and~\eqref{eq:one},
we have $\vec{B}\in C(\cD)\subseteq C(\cD)^\perp$.
Since $\vec{W}=x\in C(\cD)^\perp$ by the assumption, we have
\begin{align*}
 \vec{B}+2\vec{W}&=\vec{B\setminus W}
 \in C(\cD)^\perp\cap\{0,1\}^{33}\text{ and}\\
 \wt(\vec{B\setminus W})&=15-6=9.
\end{align*}
This contradicts the exclusion of weight $9$ established above.

\item
Suppose that $\wt(x)=12$.
Write $x=\vec{W}$, where $W\subseteq X$ and $|W|=12$.
Let $n_0,n_3,n_6,n_9$ and $n_{12}$ be as defined
in~\eqref{eq:binarymoments}.
From~\eqref{eq:distribution}, we have
\begin{align*}
 18(n_0+n_9+3n_{12})
 &=\sum_{i \in \mathcal{I}_{12}}(i-3)(i-6)n_i\\
 &=\bigl(45\cdot12+35\cdot12^2\bigr)
   -9(80\cdot12)+18\cdot176\\
 &=108.
\end{align*}
Thus, we have
\begin{equation}\label{eq:no}
 n_0+n_9+3n_{12}=6.
\end{equation}

First suppose that $n_{12}>0$.
Then there is a block $B\in\cB$ with $W\subseteq B$.
As in the case $\wt(x)=6$, we have
$\vec{B}\in C(\cD)\subseteq C(\cD)^\perp$
and $\vec{W}\in C(\cD)^\perp$.
It follows that
\begin{align*}
 \vec{B}+2\vec{W}&=\vec{B\setminus W}
 \in C(\cD)^\perp\cap\{0,1\}^{33} \text{ and}\\
 \wt(\vec{B\setminus W})&=15-12=3.
\end{align*}
This contradicts the exclusion of weight $3$ established above.
Hence, we have $n_{12}=0$.

Now suppose that $n_9>0$.
Choose a block $B\in\cB$ with $|B\cap W|=9$.
Let
\[
 u=\vec{B}+2\vec{W}\in C(\cD)^\perp.
\]
Since $|B\setminus W|=6$ and $|W\setminus B|=3$,
the vector $\phi(u)$ has $6$ entries equal to $1$,
$3$ entries equal to $-1$ and $24$ entries equal to $0$.
By Lemma~\ref{lem:signed}, we have
\[
 \langle\phi(u),\phi(\vec{B})\rangle_{\ZZ}\in\{0,3\}.
\]
On the other hand, since $u=\vec{B}+2\vec{W}\in \FF_3^{33}$,
we have $u_p=1$ for $p\in B\setminus W$
and $u_p=1+2=0$ for $p\in B\cap W$.
Thus, we have $\phi(u_p)=1$ for $p\in B\setminus W$
and $\phi(u_p)=0$ for $p\in B\cap W$.
Hence, we obtain
\[
 \langle\phi(u),\phi(\vec{B})\rangle_{\ZZ}
 =\sum_{p\in B}\phi(u_p)
 =|B\setminus W|
 =15-9=6,
\]
a contradiction.
Therefore, we have $n_9=0$.

From~\eqref{eq:no}, we have $n_0=6$.
Thus, there is a block $B\in\cB$ with $B\cap W=\emptyset$.
As in the case $\wt(x)=6$, we have
$\vec{B},\vec{W}\in C(\cD)^\perp$.
From~\eqref{eq:one}, we have $\1\in C(\cD)^\perp$.
It follows that
\begin{align*}
 \1+2\vec{B}+2\vec{W}&=\vec{X\setminus(B\cup W)}
 \in C(\cD)^\perp\cap\{0,1\}^{33} \text{ and}\\
 \wt(\vec{X\setminus(B\cup W)})&=33-15-12=6.
\end{align*}
This contradicts the exclusion of weight $6$ established above.
\end{enumerate}
This completes the proof.
\end{proof}

\subsection{Identity for a ternary code and its dual code}

For a ternary code $C$ of length $33$, define
\[
 S(C)=
 \sum_{x\in C\cap\{0,1\}^{33}}
 (-1)^{\wt(x)}(33-2\wt(x)).
\]

The following identity can be derived from the MacWilliams identity
for complete weight enumerators
(see, e.g.,~\cite[Theorem~2.3.2]{MMS})
by specialization and differentiation,
together with the assumptions $\1\in C\subseteq C^\perp$.
A direct proof of Lemma~\ref{lem:identity} using character sums
is given in Appendix~\ref{app:identity}.

\begin{lem}\label{lem:identity}
Let $C$ be a ternary code of length $33$.
Suppose that $\1\in C\subseteq C^\perp$.
Then
\begin{equation}\label{eq:identity}
 S(C^\perp)=-\frac{3^{16}}{|C|}S(C).
\end{equation}
\end{lem}

\subsection{Proof of Theorem~\ref{main}}

For a ternary code $C$ of length $33$, define the following:
\begin{align*}
 \binpart{C}&=C\cap\{0,1\}^{33}.
\end{align*}

We are ready to give a proof of Theorem~\ref{main}.

\begin{proof}[Proof of Theorem~\ref{main}]
Suppose that there exists a quasi-symmetric $2$-$(33,15,35)$ design $\cD$
with  intersection numbers $6$ and $9$.
We apply the preceding lemmas to the ternary code $C(\cD)$.

Let $E$ be either $C(\cD)$ or $C(\cD)^\perp$.
From~\eqref{eq:one},  we have $C(\cD)\subseteq C(\cD)^\perp$.
Thus, we have
$E\subseteq C(\cD)^\perp$ in either case.
If $x\in\binpart{E}$, then by Lemma~\ref{lem:weights} we have
\[
 \wt(x)\in\{0,15,18,33\}.
 \]

Let $\0$ denote the zero vector.
From~\eqref{eq:one},  we have $\1\in E$.
Thus, we have $\0,\1\in \binpart{E}$.
Hence, $\binpart{E}$ contains exactly one vector of weight $0$
and exactly one vector of weight $33$.
Therefore, the total contribution of these two vectors to $S(E)$ is
\begin{equation}\label{eq:Scounto}
 1\cdot(-1)^0(33-0)+1\cdot(-1)^{33}(33-66)
 =33+33=66.
\end{equation}
For every vector $x\in\binpart{E}$, we have $\1+2x\in E$.
Since $x$ is a $(0,1)$-vector,  we have
\[
 \1+2x\in\binpart{E}
 \text{ and }
 \wt(\1+2x)=33-\wt(x).
\]
The map $x\mapsto\1+2x$ is its own inverse, and hence is a bijection
on $\binpart{E}$.
It maps the vectors of weight $15$ to those of weight $18$,
and conversely.
Define
\[
 a_{15}(\binpart{E})
 =|\{x\in\binpart{E}\mid\wt(x)=15\}|.
\]
Therefore, $\binpart{E}$ contains exactly $ a_{15}(\binpart{E})$ vectors
of each of these two weights.

Using the definition of $S(E)$ and~\eqref{eq:Scounto}, we obtain
\begin{equation}\label{eq:Scount}
\begin{aligned}
 S(E)
&=\sum_{x\in\binpart{E}}(-1)^{\wt(x)}(33-2\wt(x))\\
 &=66+ a_{15}(\binpart{E})(-1)^{15}(33-30)
 + a_{15}(\binpart{E})(-1)^{18}(33-36)\\
 &=66-6 a_{15}(\binpart{E}).
\end{aligned}
\end{equation}
The characteristic vectors of the $176$ distinct blocks of $\cD$
are distinct codewords of weight $15$ in $C(\cD)$
whose entries are all $0$ or $1$.
Since $C(\cD)\subseteq C(\cD)^\perp$, they also belong to
$C(\cD)^\perp$.
Consequently, we have
\[
a_{15}(\binpart{C(\cD)}) \ge176 \text{ and } 
a_{15}(\binpart{C(\cD)^\perp}) \ge176.
\]
Applying~\eqref{eq:Scount} to each of the two codes, we obtain
\[
\begin{aligned}
 S(C(\cD))
 &\le66-6\cdot176=-990<0 \text{ and }\\
 S(C(\cD)^\perp)
 &\le66-6\cdot176=-990<0.
\end{aligned}
\]

On the other hand,  since \eqref{eq:one} gives $\1 \in C(\cD)\subseteq C(\cD)^\perp$,
by Lemma~\ref{lem:identity} we have
\[
 S(C(\cD)^\perp)
 =-\frac{3^{16}}{|C(\cD)|}S(C(\cD)).
\]
Since $|C(\cD)|>0$ and $S(C(\cD))<0$,
the right-hand side is positive,
contradicting $S(C(\cD)^\perp)<0$.
Therefore, no such design exists.
This completes the proof of Theorem~\ref{main}.
\end{proof}

\section{Other parameter sets}

Finally, we consider other parameter sets for exceptional
quasi-symmetric designs.

Among those listed as unresolved
in~\cite[Section~8.5.6]{BVM},
the two with the smallest numbers of points,
excluding the case settled by Theorem~\ref{main},
are those of a $2$-$(35,14,13)$ design with
intersection numbers $5$ and $8$ and a $2$-$(39,12,22)$ design
with intersection numbers $3$ and $6$.

The former is ruled out by a recent result of
Shpectorov and Zhao~\cite{SZ}.
Indeed, its block graph, in which two blocks are adjacent
when they intersect in $8$ points, would be strongly regular
with parameters $(85,14,3,2)$, contrary to their nonexistence theorem.
Together with Theorem~\ref{main}, this leaves the case of
a quasi-symmetric $2$-$(39,12,22)$ design
with intersection numbers $3$ and $6$
as the smallest remaining case.

The ternary code approach used  in this paper can also be considered
for quasi-symmetric $2$-$(v,k,\lambda)$ designs with 
intersection numbers $x$ and $y$ satisfying
\[
 k\equiv x\equiv y\equiv0\pmod3
\text{ and } r\not\equiv0\pmod3.
\]
A $2$-$(39,12,22)$ design with intersection numbers $3$ and $6$ 
satisfies these conditions, where $r=76$.
However, the present argument does not establish its nonexistence.
Indeed, suppose that such a design $\cD$ exists, and let $C(\cD)$ be its
ternary code generated by the rows of the incidence matrix of  $\cD$.
Since $ k\equiv x\equiv y\equiv0\pmod3$,
the code $C(\cD)$ is self-orthogonal.
Moreover, since $r\not\equiv0\pmod3$, $C(\cD)$ contains the all-ones vector.
Define
\[
 S_{39}(C(\cD))
 =\sum_{z\in C(\cD)\cap\{0,1\}^{39}}
 (-1)^{\wt(z)}(39-2\wt(z)),
\]
and define $S_{39}(C(\cD)^\perp)$ similarly.
The same calculation as in Appendix~\ref{app:identity} gives
\[
 S_{39}(C(\cD)^\perp)
 =(-1)^{(39+1)/2}\frac{3^{(39-1)/2}}{|C(\cD)|}S_{39}(C(\cD))
 =\frac{3^{19}}{|C(\cD)|}S_{39}(C(\cD)).
\]
The coefficient is positive, so this identity does not force
$S_{39}(C(\cD)^\perp)$ and $S_{39}(C(\cD))$ to have opposite signs.
Thus, the final sign argument used in our proof does not
yield a contradiction.
The present approach may also be useful in studying the nonexistence
of a quasi-symmetric $2$-$(45,18,34)$ design
with intersection numbers $6$ and $9$.

\bigskip
\noindent
\textbf{Acknowledgments.}
This work was supported by JSPS KAKENHI Grant Numbers 23K25784 and
25K21995.

\bigskip
\noindent
\textbf{Use of generative AI.}
The author used OpenAI's ChatGPT (GPT-6 Astra) to assist in developing
the proofs of Theorem~\ref{main} and the supporting lemmas
in November 2026.
The author independently verified all mathematical arguments
presented in the paper.

\appendix
\section{Proof of the duality identity}\label{app:identity}
In this appendix,
we give a direct proof of  Lemma~\ref{lem:identity} using character sums.

\begin{proof}[Proof of Lemma~\ref{lem:identity}]
Let $\omega\in\mathbb{C}$ be a primitive cube root of unity, so that
\[
 \omega^3=1, \omega\ne1 \text{ and }
 1+\omega+\omega^2=0.
\]
For an element $\gamma\in\FF_3$, define $\omega^\gamma=\omega^{\phi(\gamma)}$.
Since $\omega^3=1$, this value is independent of the choice
of an integer representative of $\gamma$.
For $\gamma,\delta\in\FF_3$, we have
\[
 \omega^{\gamma+\delta}=\omega^\gamma\omega^\delta,
\]
since
$\phi(\gamma+\delta)\equiv\phi(\gamma)+\phi(\delta)\pmod3$
and $\omega^3=1$.

For a vector $x\in\FF_3^{33}$, we have
\begin{equation}\label{eq:character}
 \frac{1}{|C|}\sum_{u\in C}\omega^{u\cdot x}
 =
 \begin{cases}
  1 & \text{if }x\in C^\perp,\\
  0 & \text{if }x\notin C^\perp.
 \end{cases}
\end{equation}
Indeed, if $x\in C^\perp$, then $\sum_{u\in C}\omega^{u\cdot x}=|C|$.
If $x\notin C^\perp$, choose $u_0\in C$ with $u_0\cdot x\ne0$,
and let
\[
 Q_x=\sum_{u\in C}\omega^{u\cdot x}.
\]
Since the map $u\mapsto u+u_0$ is a bijection on $C$, we have
\[
 Q_x
 =\sum_{u\in C}\omega^{(u+u_0)\cdot x}
 =\omega^{u_0\cdot x}Q_x.
\]
Since $\omega^{u_0\cdot x}\ne1$, it follows that $Q_x=0$.
This proves~\eqref{eq:character}.

For $x\in\{0,1\}^{33}$, let
\[
 f(x)=(-1)^{\wt(x)}(33-2\wt(x)).
\]
For a codeword $u\in C$, define
\[
 T(u)=\sum_{x\in\{0,1\}^{33}}f(x)\omega^{u\cdot x}.
\]
From~\eqref{eq:character} and  the definition of $S(C^\perp)$, we obtain
\begin{equation}\label{eq:transform}
\begin{aligned}
 S(C^\perp)
 &=\sum_{x\in C^\perp\cap\{0,1\}^{33}}f(x)\\
 &=\sum_{x\in C^\perp\cap\{0,1\}^{33}}f(x)\cdot1
   +\sum_{x\in\{0,1\}^{33}\setminus C^\perp}f(x)\cdot0\\
 &=\sum_{x\in\{0,1\}^{33}}
   f(x)\left(\frac{1}{|C|}
   \sum_{u\in C}\omega^{u\cdot x}\right)\\
 &=\frac{1}{|C|}\sum_{u\in C}
   \sum_{x\in\{0,1\}^{33}}f(x)\omega^{u\cdot x}\\
 &=\frac{1}{|C|}\sum_{u\in C}T(u).
\end{aligned}
\end{equation}

Fix a codeword $u=(u_1,u_2,\ldots,u_{33})\in C$.
For $j=0,1,2$, let
\[
 m_j=|\{i\in\{1,2,\ldots,33\}\mid u_i=j\}|.
\]
Since $\1\in C\subseteq C^\perp$, we have
$u\cdot u=u\cdot\1=0$.
Thus, we have
\[
 m_1+m_2\equiv0\pmod3\text{ and } m_1+2m_2\equiv0\pmod3.
\]
Since $m_0+m_1+m_2=33$, we obtain
\[
m_0\equiv0\pmod3, m_1\equiv0\pmod3 \text{ and }m_2\equiv0\pmod3.
\]

Consider the following polynomial:
\[
 H_u(z)
 =\prod_{p=1}^{33}(1+\omega^{u_p} z)
 =(1+z)^{m_0}(1+\omega z)^{m_1}(1+\omega^2z)^{m_2}.
\]
For each $x=(x_1,x_2,\ldots,x_{33})\in\{0,1\}^{33}$,
choose $1$ from the factor $1+z\omega^{u_p}$ if $x_p=0$
and choose $z\omega^{u_p}$ if $x_p=1$.
The resulting product is
\[
 \prod_{\substack{1\le p\le33\\x_p=1}}z\omega^{u_p}
 =z^{\wt(x)}\omega^{\sum_{p=1}^{33}u_px_p}
 =z^{\wt(x)}\omega^{u\cdot x}.
\]
These choices are in bijection with the vectors
$x\in\{0,1\}^{33}$.
Thus, summing over all such choices, we obtain
\[
 H_u(z)
 =\sum_{x\in\{0,1\}^{33}}z^{\wt(x)}\omega^{u\cdot x}.
\]

Differentiating $H_u(z)$ and multiplying by $z$, we obtain
\[
 zH'_u(z)
 =\sum_{x\in\{0,1\}^{33}}
   \wt(x)z^{\wt(x)}\omega^{u\cdot x}.
\]
Hence, we have
\[
 33H_u(z)-2zH'_u(z)
 =\sum_{x\in\{0,1\}^{33}}
   \bigl(33-2\wt(x)\bigr)z^{\wt(x)}\omega^{u\cdot x}.
\]
Therefore, setting $z=-1$ and using the definitions of $f$ and $T(u)$,
we obtain
\begin{equation}\label{eq:derivative}
\begin{aligned}
 T(u)
 &=\sum_{x\in\{0,1\}^{33}}f(x)\omega^{u\cdot x}\\
 &=\sum_{x\in\{0,1\}^{33}}
   (-1)^{\wt(x)}\bigl(33-2\wt(x)\bigr)\omega^{u\cdot x}\\
 &=33H_u(-1)+2H'_u(-1).
\end{aligned}
\end{equation}

We consider the two cases $m_0>0$ and $m_0=0$ separately.
\begin{enumerate}
\item
Suppose that $m_0>0$.
Since $m_0$ is divisible by $3$, we have $m_0\ge3$.
Thus, $H_u(z)$ is divisible by $(1+z)^3$, and consequently we have
\[
 H_u(-1)=H'_u(-1)=0.
\]
From~\eqref{eq:derivative}, we have $T(u)=0$.

\item
Suppose that $m_0=0$.
Let
\[
 \alpha=1-\omega\text{ and } \beta=1-\omega^2.
\]
From $1+\omega+\omega^2=0$ and $\omega^3=1$, we obtain
\begin{equation}\label{eq:alpha}
 \beta=-\omega^2\alpha,
 \alpha^2=-3\omega \text{ and }
 \alpha^{32}=(-3\omega)^{16}=3^{16}\omega.
\end{equation}
Since $m_2$ is divisible by $3$, we have $\omega^{2m_2}=1$.
Since $m_0=0$, we have $m_1+m_2=33$.
Substituting $z=-1$ into $H_u(z)$ and using~\eqref{eq:alpha}
together with $\omega^{2m_2}=1$, we obtain
\begin{equation}\label{eq:Hvalue}
\begin{aligned}
 H_u(-1)&=(1-\omega)^{m_1}(1-\omega^2)^{m_2}=\alpha^{m_1}\beta^{m_2}
 =\alpha^{m_1+m_2}(\alpha^{-1}\beta)^{m_2}\\
 & =\alpha^{m_1+m_2}(-\omega^2)^{m_2}
=(-1)^{m_2}\alpha^{33}.
\end{aligned}
\end{equation}
Since $m_0=0$, every coordinate of $u$ is either $1$ or $2$.
Thus, none of the factors of $H_u(z)$ vanishes at $z=-1$.
Hence, $H_u(z)\ne0$ in a neighbourhood of $-1$.
By the product rule, we have
\[
 H'_u(z)
 =\sum_{p=1}^{33}\omega^{u_p}
   \prod_{\substack{1\le q\le33\\q\ne p}}
   (1+\omega^{u_q}z).
\]
Dividing both sides by $H_u(z)$ and cancelling common factors,
we obtain
\[
 \frac{H'_u(z)}{H_u(z)}
 =\sum_{p=1}^{33}
   \frac{
     \omega^{u_p}
     \displaystyle\prod_{\substack{1\le q\le33\\q\ne p}}
     (1+\omega^{u_q}z)}
   {\displaystyle\prod_{q=1}^{33}(1+\omega^{u_q}z)}
 =\sum_{p=1}^{33}
   \frac{\omega^{u_p}}{1+\omega^{u_p}z}.
\]
Evaluating at $z=-1$ and using~\eqref{eq:alpha}
and~\eqref{eq:Hvalue}, we obtain
\begin{equation}\label{eq:Hderivative}
\begin{aligned}
 H'_u(-1)
 &=H_u(-1)
   \left(\frac{m_1\omega}{\alpha}
         +\frac{m_2\omega^2}{\beta}\right)\\
 &=(-1)^{m_2}\alpha^{32}(m_1\omega-m_2).
\end{aligned}
\end{equation}

From~\eqref{eq:derivative}, \eqref{eq:alpha},
\eqref{eq:Hvalue} and~\eqref{eq:Hderivative},
we obtain
\begin{equation}\label{eq:Tu}
\begin{aligned}
T(u)
&=(-1)^{m_2}\alpha^{32}
  (33\alpha+2m_1\omega-2m_2)\\
 &=(-1)^{m_2}3^{16}
   (33\omega\alpha+2m_1\omega^2-2m_2\omega)\\
 &=(-1)^{m_2}3^{16}
   ((33-2m_2)\omega+(2m_1-33)\omega^2)\\
 &=(-1)^{m_2}3^{16}(m_1-m_2)(\omega+\omega^2)\\
 &=-3^{16}(-1)^{m_2}(33-2m_2).
 \end{aligned}
\end{equation}
Here, the fourth equality follows from $m_1+m_2=33$,
and the last equality follows from
$\omega+\omega^2=-1$ and $m_1-m_2=33-2m_2$.

Since $m_0=0$, the vector $u+2\1$ has $m_1$ entries equal to $0$
and $m_2$ entries equal to $1$.
Thus, $u+2\1$ is a $(0,1)$-vector of weight $m_2$.
Hence, the definition of $f$ gives
\[
 f(u+2\1)=(-1)^{m_2}(33-2m_2).
\]
From~\eqref{eq:Tu}, we obtain
\[
 T(u)=-3^{16}f(u+2\1).
\]
\end{enumerate}

Combining the two cases $m_0>0$ and $m_0=0$, we have
\begin{equation}\label{eq:Tvalue}
 T(u)=
 \begin{cases}
  0 & \text{if }\wt(u)<33,\\
  -3^{16}f(u+2\1) & \text{if }\wt(u)=33.
 \end{cases}
\end{equation}
Since $\1\in C$, the map $u\mapsto u+2\1$ sends every codeword
of weight $33$ in $C$ to a vector of $\binpart{C}$.
Conversely, if $x\in\binpart{C}$, then $u=x+\1$ belongs to $C$
and has every coordinate equal to $1$ or $2$.
Thus, we have $\wt(u)=33$.
These maps are inverse to each other.
Hence, $u\mapsto u+2\1$ is a bijection from the set of codewords
of weight $33$ in $C$ onto $\binpart{C}$.
Using this bijection and~\eqref{eq:Tvalue}, we obtain
\[
\begin{aligned}
 \sum_{u\in C}T(u)
 &=\sum_{\substack{u\in C\\\wt(u)<33}}T(u)
   +\sum_{\substack{u\in C\\\wt(u)=33}}T(u)\\
 &=0-3^{16}
   \sum_{\substack{u\in C\\\wt(u)=33}}f(u+2\1)\\
 &=-3^{16}\sum_{x\in\binpart{C}}f(x)\\
 &=-3^{16}S(C).
\end{aligned}
\]
From~\eqref{eq:transform}, we obtain
\[
 S(C^\perp)
 =\frac{1}{|C|}\sum_{u\in C}T(u)
 =-\frac{3^{16}}{|C|}S(C),
\]
as required.
\end{proof}



\end{document}